\documentclass{amsart}

\usepackage{graphicx,amssymb,mathrsfs,amsmath,color,fancyhdr}
\usepackage[all]{xy}
\newtheorem{theorem}{Theorem}[section]
\newtheorem{lemma}[theorem]{Lemma}

\theoremstyle{definition}

\newtheorem{example}[theorem]{Example}

\newtheorem{proposition}[theorem]{Proposition}

\theoremstyle{remark}

\numberwithin{equation}{section}

\begin{document}

\title[On $*$-clean group rings over finite fields]
{On $*$-clean group rings over finite fields}

\author{Dongchun Han}
\address{Department of Mathematics, Southwest Jiaotong University, Chengdu 610000, P.R. China}
\email{handongchun@swjtu.edu.cn}
\author{Hanbin Zhang}
\address{School of Mathematics (Zhuhai), Sun Yat-sen University, Zhuhai 519082,
Guangdong, P.R. China}
\email{zhanghb68@mail.sysu.edu.cn}

\begin{abstract}
A ring $R$ is called {\sl clean} if every element of $R$ is the sum of a unit and an idempotent. Motivated by a question proposed by Lam on the cleanness of von Neumann Algebras, Va\v{s} introduced a more natural concept of cleanness for $*$-rings, called the $*$-cleanness. More precisely, a $*$-ring $R$ is called a $*$-$clean$ $ring$ if every element of $R$ is the sum of a unit and a projection ($*$-invariant idempotent). Let $\mathbb F$ be a finite field and $G$ a finite abelian group. In this paper, we introduce two classes of involutions on group rings of the form $\mathbb FG$ and characterize the $*$-cleanness of these group rings in each case. When $*$ is taken as the classical involution, we also characterize the $*$-cleanness of $\mathbb F_qG$ in terms of LCD abelian codes and self-orthogonal abelian codes in $\mathbb F_qG$.
\end{abstract}

\keywords{group ring; *-cleanness; Galois theory; primitive idempotent; abelian group code}

\maketitle{}

\section{Introduction}

A ring $R$ is called {\sl clean} if every element of $R$ is the sum of a unit and an idempotent. In 1977, Nicholson \cite{[N1]} introduced the clean rings and related them to exchange rings. A clean ring can be regarded as an additive analog of a unit-regular ring in which each element is the product of a unit and an idempotent. Some important examples of clean rings include local rings, semiperfect rings and left (right) Artinian rings. Many interesting results have been established about clean rings; see, e.g., \cite{CNZ,CWZ,LZ,LZ1,Mc1,Zhou}.

In 2005, T.Y. Lam \cite[Introduction]{[V]} asked which von Neumann algebras are clean as rings; see \cite{KR} for a detailed introduction of von Neumann algebras and \cite{CHWYZ} for a recent progress. Trying to answer this question, Va\v{s} \cite{[V]} pointed out that it is more natural to utilize the fact that a von Neumann algebra is a $*$-ring and that the projections are $*$-invariant idempotents. In fact, a ring $R$ is called a $*$-$ring$ (or ring with involution $*$) if there exists an operation $*:R\rightarrow R$ such that
$$(x+y)^*=x^*+y^*,\text{ }(xy)^*=y^*x^*,\text{ and }(x^*)^*=x,$$
for all $x,y\in R$. We call an element $p$ of a $*$-ring $R$ a $projection$ (motivated by the definition of the projections in operator algebras) if $p$ is a $*$-invariant idempotent, i.e., $p^*=p=p^2$, and call a $*$-ring $R$ a $*$-$clean$ $ring$ if each element of $R$ is the sum of a unit and a projection. In \cite{[V]}, Va\v{s} also proposed a question that whether there exists a clean ring (with involution $*$) that is not $*$-clean. A year later, Li and Zhou {\rm\cite{[lz]}} answered this question affirmatively.

Meanwhile, the study of the $*$-cleanness of group rings has received a lot of attention. Let $R$ be a ring and $G$ a group. We denote by $RG$ the group ring of $G$ over $R$. It is well-known that for a commutative ring $R$, the map $*:RG\rightarrow RG$ given by $(\sum a_g g)^*=\sum a_gg^{-1}$ is an involution which is called the classical (or standard) involution on $RG$; see \cite[Proposition 3.2.11]{[MS]}. In {\rm\cite{[GCL]}}, Gao, Chen and Li characterized the $*$-cleanness of group rings $RG$, where $R$ is a commutative local ring and $G$ is of small orders. Later, Li, Parmenter and Yuan {\rm\cite{[lpy]}} and Huang, Li and Yuan {\rm\cite{[HLY]}} studied when a group ring $\mathbb{F}C_n$ is $*$-clean, where $\mathbb{F}$ is a field and $C_n$ is the cyclic group of order $n$ (also see \cite{HLT} for a study in the non-commutative case). Later, Han, Ren and Zhang \cite{HRZ} and Tang, Wu and Li \cite{TWL} independently extended the above results to the case $\mathbb{F}_qG$, where $\mathbb{F}_q$ is a finite field with $q$ elements and $G$ is a finite abelian group with exponent $n$ and $(q,n)=1$. More precisely, in \cite{HRZ}, it is proved that $\mathbb F_qG$ is $*$-clean if and only if there exists $\sigma\in Gal(\mathbb F_q(\omega_n)/\mathbb F_q)$ such that $\sigma(\omega_n) = \omega_n^{-1}$, where $\omega_n$ is an $n$-th primitive root of unity over $\mathbb F_q$. Equivalently, $\mathbb F_qG$ is $*$-clean if and only if there exists a positive integer $t$ such that $q^t\equiv -1\pmod{n}$; see Section 3 for a discussion. Recently, Wu, Yue and Tang \cite{WYT} introduced another involution $*:\ \mathbb F_{q^2}G\rightarrow \mathbb F_{q^2}G$, $\sum_{g\in G}a_g g\mapsto\sum_{g\in G}a_g^q g^{-1}$ and characterized the $*$-cleanness of $\mathbb F_{q^2}G$. We also refer to \cite{CH} for a further study of the above characterization.

In this paper, using a different approach, we shall study two new classes of involutions and extend the previous results. Let $\mathbb F_q$ be a finite field with $q$ elements of characteristic $p$. Let $G$ be a finite abelian group with exponent $n$. The first class of involutions we are going to study is a generalized version of the classical one. Let $\sigma_1: \mathbb F_{q}G\rightarrow \mathbb F_{q}G$ with $\sum_{g\in G}a_gg\mapsto\sum_{g\in G}a_gg^v$, where $v^2\equiv 1\pmod{n}$ and $v\not\equiv 1\pmod{n}$. Then we have the following characterization for this class of involutions.
\begin{theorem}\label{mainthm}
Let $\mathbb F_q$ be a finite field with $q$ elements of characteristic $p$ and $G$ a finite abelian group with exponent $n$. Assume that $G=Syl_p(G)\times H$ and $m$ is the exponent of $H$, where $Syl_p(G)$ is the Sylow $p$-subgroup of $G$. Taking the involution $\sigma_1$ as $*$ on $\mathbb{F}_{q}G$, then $\mathbb{F}_qG$ is $*$-clean if and only if there exists a positive integer $t$ such that $q^t\equiv v\pmod{m}$.
\end{theorem}

In Proposition \ref{all involutions}, we shall show that there are infinitely many group rings $\mathbb{F}_qG$ on which the only possible involutions are of the form $\sigma_1$. Moreover, we shall provide an elementary method to study the $*$-cleanness of $\mathbb F_2G$.

Let $\mathbb F_{q^2}$ be a finite field with $q^2$ elements. The second class of involutions we are going to study is the following: Let $\sigma_2: \mathbb F_{q^2}G\rightarrow \mathbb F_{q^2}G$ with $\sum_{g\in G}a_gg\mapsto\sum_{g\in G}a_g^{q}g^v$, where $v^2\equiv 1\pmod{n}$. Then we have the following characterization for this class of involutions.

\begin{theorem}\label{mainthm2}
Let $\mathbb F_{q^2}$ be a finite field with $q^2$ elements of characteristic $p$ and $G$ a finite abelian group with exponent $n$. Assume that $G=Syl_p(G)\times H$ and $m$ is the exponent of $H$, where $Syl_p(G)$ is the Sylow $p$-subgroup of $G$. Taking the involution $\sigma_2$ as $*$ on $\mathbb{F}_{q^2}G$, then $\mathbb{F}_{q^2}G$ is $*$-clean if and only if there exists a positive integer $t$ such that $q^{2t}\equiv qv\pmod{m}$.
\end{theorem}

For the classical involution $*$, an important application of the study of the $*$-cleanness of $\mathbb F_qG$ can be found in coding theory. We briefly recall some basic definitions in coding theory. An abelian code over $\mathbb F_q$ is defined as an ideal in the group ring $\mathbb F_qG$. Note that when $G$ is a cyclic group, then these codes are the well-known and widely used cyclic codes. The study of abelian codes or, more generally, group codes (which are one sided ideals in arbitrary finite group algebras) has received a lot of attention; see, e.g., \cite{KelS, MacW}. We refer to \cite{Ward} for a discussion on non-cyclic group codes with their applications. For $a=\sum_{g\in G}a_g g, b=\sum_{g\in G}b_g g\in \mathbb F_qG$, we define the following inner product on $\mathbb F_qG$: $\langle a,b\rangle
=\sum_{g\in G}a_gb_g$. A linear code $\mathcal C$ is called an LCD code  (linear code with complementary dual) if $\mathcal C\cap \mathcal C^{\bot}=\{0\}$, where
$$\mathcal C^{\bot}=\{b\in\mathbb F_qG\ |\ \langle a,b\rangle=0\text{ for all }a\in\mathcal C\}.$$
If $\mathcal C\subset \mathcal C^{\bot}$, then $\mathcal C$ is called a self-orthogonal code. LCD codes have been widely applied in data storage, communications systems, consumer electronics, and cryptography. Self-orthogonal codes also have wide applications in communications and secret sharing. For more background and some recent progress on LCD and self-orthogonal code, we refer to, e.g., \cite{LYW,ZLTD} and the references therein. In Section 4, when $*$ is taken as the classical involution, we shall characterize the $*$-cleanness of $\mathbb F_qG$ in terms of LCD abelian codes and self-orthogonal abelian codes in $\mathbb F_qG$.

The following sections are organized as follows. In Section 2, we shall introduce some definitions and auxiliary results. In Section 3, we prove Theorems \ref{mainthm} and \ref{mainthm2}. Moreover, we provide an elementary method to study the $*$-cleanness of $\mathbb F_2G$. In Section 4, we discuss the connection between the $*$-cleanness of $\mathbb F_qG$ and abelian group codes in $\mathbb F_qG$.

\section{Preliminaries}

This section will provide more rigorous definitions and some notation. We also introduce some auxiliary results that will be used repeatedly below.

Let $\mathbb{N}$ denote the set of all positive integers. Throughout this paper, let $G$ be an abelian group written multiplicatively. By the fundamental theorem of finite abelian groups we have
$$G\cong C_{m_1}\times\cdots\times C_{m_r}\cong\langle x_1\rangle\times\langle x_2\rangle\times\cdots\times\langle x_r\rangle,$$
where $m_1,\ldots,m_r\in\mathbb{N}$ with $1<m_1\mid\ldots\mid m_r$. Moreover, $m_1,\ldots,m_r$ are uniquely determined by $G$, and $m_r$ is called the $exponent$ of $G$.

Let $R$ be a unitary ring (with the identity $1_R$) and $G$ a group. We denote by $RG$ the group ring of $G$ over $R$. We regard $R$ as a subring of $RG$ via the natural isomorphism $R\cong R1_G$, where $1_G$ is the identity of $G$. Let $\mathbb{F}$ be a field. By the basic Galois theory, $\mathbb{F}(\omega)$ is a Galois extension of $\mathbb{F}$, where $\omega$ is an $n$-th primitive root of unity over $\mathbb F$. For an element $a=\sum_{g\in G}a_gg\in \mathbb{F}(\omega)G$, we define that $\sigma(a)=\sum_{g\in G}\sigma(a_g)g$, where $\sigma\in Gal(\mathbb{F}(\omega)/\mathbb{F})$. Note that for $a,b\in\mathbb{F}(\omega)G$, we have $\sigma(ab)=\sigma(a)\sigma(b)$ for any $\sigma\in Gal(\mathbb{F}(\omega)/\mathbb{F})$.

Recall that, in a group ring, an idempotent is called primitive if it can not be written as the sum of two non-zero idempotents.

\begin{lemma} \label{lemma3.1}{\rm(\cite{[lz]}, Theorem 2.2)}
A commutative $*$-ring is $*$-clean if and only if it is clean and every idempotent is a projection.
\end{lemma}

The proof of the following lemma is routine, which we omit here.

\begin{lemma}\label{sumRG}
Let $\{R_i\}_{i\in I}$ be a family of rings and set $R=\oplus_{i\in I} R_i$. Then for any group $G$, we have that
$$RG\cong\oplus_{i\in I}R_iG.$$
\end{lemma}

\begin{lemma}\label{cleanFG}
Let $\mathbb{F}$ be a field of characteristic $p>0$ and $G$ a finite abelian group of order $p^km$ with $(p,m)=1$. Let $Syl_{p}(G)$ be the Sylow $p$-subgroup of $G$ and $H$ the subgroup of order $m$. Then $\mathbb FG$ is a clean ring.
\end{lemma}

\begin{proof}
We assume that $\mathbb FH\cong\mathbb K^1\oplus\cdots\oplus\mathbb K^s$, where the $\mathbb K^i$'s are finite extensions of $\mathbb F$. Moreover, by Lemma \ref{sumRG} we have
$$\mathbb FG=\mathbb F(H\times Syl_{p}(G))\cong(\mathbb FH)Syl_p(G)\cong
\mathbb K^1Syl_p(G)\oplus\cdots\oplus\mathbb K^sSyl_p(G).$$
It suffices to prove that each $\mathbb K^iSyl_p(G)$ is a clean ring.
Let $\alpha=\sum_{g\in Syl_p(G)}a_gg\in \mathbb K^iSyl_p(G)$. We claim that $\alpha$ is invertible if and only if $\sum_{g\in Syl_p(G)}a_g^{p^{k}}\neq 0$. If $\alpha$ is invertible, then $\alpha^{p^{k}}=\sum_{g\in Syl_p(G)}a_g^{p^{k}}$ is also invertible. Conversely, if $\sum_{g\in Syl_p(G)}a_g^{p^{k}}\neq 0$, then $\alpha^{p^{k}}$ is invertible. Therefore, there exists $b\in \mathbb K^iSyl_p(G)$ such that $\alpha^{p^{k}}b=1$. Consequently, $\alpha^{p^{k}-1}b$ is the inverse of $\alpha$. If $\alpha$ is a unit, then $\alpha=\alpha+0$ is a representation of $\alpha$ as the sum of a unit and an idempotent. Otherwise, by the above claim, $\alpha-1$ is a unit and $\alpha=(\alpha-1)+1$ is a representation of $\alpha$ as the sum of a unit and an idempotent. This completes the proof.
\end{proof}

Note that, when $(p,|G|)=1$, the cleanness of $\mathbb FG$ follows immediately from the fact that it is semi-simple.

For any two positive integers $n,m$ with $(n,m)=1$, we denote by ord$_n(m)$ the order of $m$ in the multiplicative group $(\mathbb Z/n\mathbb Z)^{\times}$. We shall use the following lemma from elementary number theory, a proof can be found in \cite[Theorem 3.6]{[Nathanson]}.

\begin{lemma}\label{ord_p}
Let $a,n\in\mathbb N$ and $a\neq\pm 1$. Let $p$ be an odd prime, $d=\text{ord}_p(a)$ and $p^h||a^d-1$. Then
    \[
 ord_{p^n}(a)= \left\{ \begin{array}{ll} & d,\quad\quad \mbox{ if } 1\le n\le h; \\ &
 p^{n-h}d, \mbox{ if } n\ge h.
\end{array} \right.
\]
\end{lemma}

\section{Proofs of the main results}

In this section, we shall prove Theorems \ref{mainthm} and \ref{mainthm2}. In \cite[Section 5]{HRZ}, although somewhat complicated, an explicit construction of the primitive idempotents in $\mathbb FG$ was given, where $\mathbb{F}$ is a field of characteristic $p>0$ and $G$ is a finite abelian group with exponent $n$ and $(p,n)=1$. In the following, we provide a more general and simplified construction.

Let $\mathbb{F}$ be a field of characteristic $p>0$ and $G$ a finite abelian group of order $p^km$ with $(p,m)=1$. Let $Syl_{p}(G)$ be the Sylow $p$-subgroup of $G$ and $H$ the subgroup of order $m$. Then clearly we have $G=Syl_{p}(G)\times H$. In the following lemma, we show that the problem of finding all the primitive idempotents in $\mathbb FG$ is reduced to the same problem in $\mathbb FH$.

\begin{lemma}\label{FGFH}
Let $\mathbb{F}$ be a field of characteristic $p>0$ and $G$ a finite abelian group of order $p^km$ with $(p,m)=1$. Assume that $G=Syl_p(G)\times H$, where $Syl_p(G)$ is the Sylow $p$-subgroup of $G$ and $H$ is the subgroup of order $m$. Then we have any idempotent $e$ in $\mathbb FG$ actually belongs to $\mathbb FH$. In other words, $\mathbb FH$ contains exactly all idempotents in $\mathbb FG$.
\end{lemma}

\begin{proof}
Since $e^2=e$, it follows that $e^i=e$ for any integer $i\ge 2$. In particular, $e^{p^k}=e$. Let $e=\sum_{g\in G}a_gg$. As $e^{p^k}=(\sum_{g\in G}a_gg)^{p^k}=\sum_{g\in G}a_g^{p^k}g^{p^k}\in\mathbb FH$, the desired result follows.
\end{proof}

With Lemma \ref{FGFH}, it suffices to construct all the primitive idempotents in $\mathbb FH$, where $\mathbb{F}$ is a field of characteristic $p>0$ and $H$ is a finite abelian group of order $m$ with $(p,m)=1$. The construction below is similar to that in \cite{BrRio}, we provide it here for the convenience of readers.

Let $\overline{\mathbb F}$ be an algebraic closure of $\mathbb F$. We consider the group $\hat{H}$ which consists of all characters of $H$ over $\overline{\mathbb F}$, that is,
$$\hat{H}=\{\psi\ |\ \psi:H\rightarrow \overline{\mathbb F}\text{ a homomorphism}\}.$$
It is well-known that $|\hat{H}|=|H|$. For each $\psi\in\hat{H}$, we define
\begin{equation}\label{psi}
e_{\psi}=\frac{1}{|H|}\sum_{h\in H}\psi(h)h,
\end{equation}
which is an element in $\overline{\mathbb F}H$. Clearly, we have the following three properties:
\begin{enumerate}
\item $e_{\psi}^2=e_{\psi}$, for any $\psi\in \hat{H}$;
\item $e_{\psi}e_{\varphi}=0$, for any $\psi, \varphi\in \hat{H}$ with $\psi\neq\varphi$;
\item $\sum_{\psi\in\hat{H}}e_{\psi}=1$.
\end{enumerate}
Therefore, the set
$$\mathcal{E}=\{e_\psi\ |\ \psi\in\hat{H}\}$$ contains exactly all primitive idempotents of $\overline{\mathbb F}H$. Then we can construct the primitive idempotents of $\mathbb FH$ from the set $\mathcal{E}$. For any fixed $\psi\in \hat{H}$, let $d$ be the order of $\psi$ in $\hat{H}$ and $\omega_d$ a $d$-th primitive root of unity over $\mathbb F$. Then we have $e_{\psi}\in \mathbb F(\omega_d)H$. For simplicity, we define
\begin{equation}\label{Trace}
\text{Tr}_{\psi}(e_\psi):=\text{Tr}_{\mathbb F(\omega_d)/\mathbb F}(e_{\psi}).
\end{equation}
We have the following lemma.
\begin{lemma}\label{AllPIFH}
Let $\psi\in\hat{H}$. Then we have $\text{Tr}_{\psi}(e_\psi)$ is a primitive idempotent in $\mathbb FH$. Moreover, the set $$\mathcal{F}:=\{\text{Tr}_{\psi}(e_\psi)\ |\ \psi\in\hat{H}\}$$ contains exactly all primitive idempotents of $\mathbb FH$.
\end{lemma}
\begin{proof}
It is easy to see that $\sigma(e_\psi)$ is also a primitive idempotent in $\mathbb F(\omega_d)H$, for any $\sigma\in Gal(\mathbb F(\omega_d)/\mathbb F)$. Moreover, we have $\sigma(e_\psi)\neq\sigma'(e_\psi)$ for any $\sigma\neq\sigma'$. Therefore, due to the properties (1)-(3) of $\mathcal E$, $\text{Tr}_{\psi}(e_\psi)$ is an idempotent in $\mathbb FH$. For any idempotent $u$ in $\mathbb FH$, we have $u=e_{\psi_1}+\cdots+e_{\psi_k}$. Clearly, $\sigma(e_{\psi_1})$ also should appear in the summation of $u$ for any $\sigma\in Gal(\mathbb F(\omega_d)/\mathbb F)$. Consequently, $\text{Tr}_{\psi}(e_\psi)$ is a primitive idempotent in $\mathbb FH$. The second part of this lemma follows from the above discussion.
\end{proof}

By Lemmas \ref{FGFH} and \ref{AllPIFH}, we have the following.

\begin{proposition}\label{PIFG}
Let $\mathbb{F}$ be a field of characteristic $p>0$ and $G$ a finite abelian group of order $p^km$ with $(p,m)=1$. Let $H$ be the subgroup of $G$ of order $m$. Then the set $$\mathcal{F}:=\{\text{Tr}_{\psi}(e_\psi)\ |\ \psi\in\hat{H}\}$$ contains exactly all primitive idempotents of $\mathbb FG$.
\end{proposition}

Note that the above mentioned $\mathbb F$ is not necessarily finite. In the following, we consider the case when $\mathbb F=\mathbb F_q$, a finite field with $q$ elements. We have the following refined description of the set $\mathcal F$ mentioned in the above lemma.

\begin{lemma}\label{PIC}
Let $\mathbb F=\mathbb F_q$ be a finite field with $q$ elements. Let $\psi_1,\psi_2\in\hat{H}$. Then $\text{Tr}_{\psi_1}(e_{\psi_1})= \text{Tr}_{\psi_2}(e_{\psi_2})$ if and only if $\psi_1=\psi_2^{q^t}$ holds for some $t\in\mathbb N$.
\end{lemma}

\begin{proof}
We assume that the order of $\psi_i$ is $d_i$, for $i=1,2$. Then
$\text{Tr}_{\psi_1}(e_{\psi_1})=\sum_{\sigma\in Gal(\mathbb F_q(\omega_{d_1})/\mathbb F_q)}\sigma(e_{\psi_1})$ and $\text{Tr}_{\psi_2}(e_{\psi_2})=\sum_{\tau\in Gal(\mathbb F_q(\omega_{d_2})/\mathbb F_q)}\tau(e_{\psi_2})$. Since for any $\sigma\in Gal(\mathbb F_q(\omega_{d_1})/\mathbb F_q)$ and $\tau\in Gal(\mathbb F_q(\omega_{d_2})/\mathbb F_q)$, $\sigma(e_{\psi_1})$ and $\tau(e_{\psi_2})$ are primitive idempotents in $\overline{\mathbb F_q}G$. If $\text{Tr}_{\psi_1}(e_{\psi_1})= \text{Tr}_{\psi_2}(e_{\psi_2})$, then we must have $e_{\psi_1}=\tau(e_{\psi_2})$ holds for some $\tau\in Gal(\mathbb F_q(\omega_{d_2})/\mathbb F_q)$. Therefore,
\begin{equation}\label{psitau}
\psi_1(g)=\tau\psi_2(g)
\end{equation}
holds for any $g\in H$. Moreover, as $Gal(\mathbb F_q(\omega_{d_2})/\mathbb F_q)$ is generated by the Frobenius map: $$\Phi:\mathbb{F}_q(\omega_{d_2})\rightarrow \mathbb{F}_q(\omega_{d_2}),\text{ and }\Phi(x)=x^q,$$ together with (\ref{psitau}) we have that $\psi_1=\psi_2^{q^t}$ for some $t\in\mathbb N$. Conversely, we assume that $\psi_1=\psi_2^{q^t}$ holds for some $t\in\mathbb N$. Since $\text{Tr}_{\psi_2}(e_{\psi_2})\in\mathbb F_qH$, we have $\text{Tr}_{\psi_2}(e_{\psi_2})=\text{Tr}_{\psi_2^{q^i}}(e_{\psi_2^{q^i}})$ holds for any $i\in\mathbb N$. In particular, we have $\text{Tr}_{\psi_1}(e_{\psi_1})=\text{Tr}_{\psi_2}(e_{\psi_2})$. This completes the proof.
\end{proof}

Let $d$ be the order $\psi\in\hat{H}$ and $k$ the multiplicative order of $q$ modulo $d$. The set $\{\psi,\psi^{q},\ldots,\psi^{q^{k-1}}\}$ is usually called a $q$-cyclotomic class in $\hat{H}$ in the literature; see, e.g., \cite{BrRio}. Therefore, the $q$-cyclotomic classes in $\hat{H}$ correspond to all primitive idempotents in $\mathbb FH$. Now, we can prove Theorem \ref{mainthm} using the above construction.

\medskip

{\sl Proof of Theorem \ref{mainthm}.}
Firstly, we assume that $\mathbb{F}_qG$ is $*$-clean.  Let $\psi$ be an element of order $m$ in $\hat{H}$.
By Lemmas \ref{lemma3.1} and \ref{FGFH}, we have $u=u^*$ for all idempotents $u\in \mathbb{F}_qH$, and in particular $\text{Tr}_{\psi}(e_\psi)={\text{Tr}_{\psi}(e_\psi)}^*$. Moreover, we have
\begin{align*}
{\text{Tr}_{\psi}(e_\psi)}^*&=\sum_{\sigma\in Gal(\mathbb F_q(\omega_{n})/\mathbb F_q)}\sigma\left(\frac{1}{|H|}\sum_{h\in H}\psi(h)h\right)^*\\
&=\sum_{\sigma\in Gal(\mathbb F_q(\omega_{n})/\mathbb F_q)}\sigma\left(\frac{1}{|H|}\sum_{h\in H}\psi(h)h^v\right)\\
&=\sum_{\sigma\in Gal(\mathbb F_q(\omega_{n})/\mathbb F_q)}\sigma\left(\frac{1}{|H|}\sum_{h\in H}\psi(h^{v^2})h^v\right)\\
&=\sum_{\sigma\in Gal(\mathbb F_q(\omega_{n})/\mathbb F_q)}\sigma\left(\frac{1}{|H|}\sum_{h\in H}\psi^v(h)h\right)=\text{Tr}_{{\psi^v}}(e_{\psi^v}).
\end{align*}
Therefore $\text{Tr}_{\psi}(e_\psi)=\text{Tr}_{{\psi^v}}(e_{\psi^v})$, and by Lemma \ref{PIC} we have $\psi^v=\psi^{q^t}$ for some $t\in\mathbb N$. Based on our choice of $\psi$, the desired result follows.

Conversely, we assume that there exists $t\in\mathbb N$ such that $q^t\equiv v\pmod{m}$. By Lemma \ref{cleanFG}, $\mathbb{F}_qG$ is clean. By Lemma \ref{lemma3.1}, we only need to prove that every idempotent is a projection. Let $u$ be any idempotent in $\mathbb{F}_qG$, then $u=\text{Tr}_{\psi_{i_1}}(e_{\psi_{i_1}})+\cdots+
\text{Tr}_{\psi_{i_t}}(e_{\psi_{i_l}})$, for some $\psi_{i_1},\ldots,\psi_{i_l}\in\hat{H}$.
Therefore, it suffices to prove that $\text{Tr}_{\psi}(e_{\psi})=\text{Tr}_{\psi}(e_{\psi})^*$ holds for any $\psi\in\hat{H}$. For any fixed $\psi\in\hat{H}$, let $d$ be the order of $\psi$ in $\hat{H}$ and $\omega_d$ a $d$-th primitive root of unity over $\mathbb F_q$. Since there exists $t\in\mathbb N$ such that $q^t\equiv v\pmod{m}$, we have $q^t\equiv v\pmod{d}$. This is equivalent to say that there exists $\sigma\in Gal(\mathbb{F}_q(\omega_{d})/\mathbb{F}_q)$ such that $\sigma(\omega_{d})=\omega_{d}^{v}$. We assume that $Gal(\mathbb{F}_q(\omega_{d})/\mathbb{F}_q)=\{\theta_1,\ldots,\theta_{t}\}$. Therefore we have
\begin{align*}
\text{Tr}_{\psi}(e_{\psi})&=\text{Tr}_{\mathbb F_q(\omega_d)/\mathbb F_q}(e_{\psi})=\sum_{i=1}^{t}\theta_i(e_{\psi})
=\sigma\left(\sum_{i=1}^{t}\theta_i(e_{\psi})\right)\\
&=\frac{1}{|H|}\left(\sum_{i=1}^{t}\sigma\theta_i\big(\sum_{h\in H}\psi(h)h\big)\right)
=\frac{1}{|H|}\left(\sum_{i=1}^{t}\theta_i\sigma\big(\sum_{h\in H}\psi(h)h\big)\right)\\
&=\frac{1}{|H|}\left(\sum_{i=1}^{t}\theta_i\big(\sum_{h\in H}\psi(h)^vh\big)\right)=
\frac{1}{|H|}\left(\sum_{i=1}^{t}\theta_i\big(\sum_{h\in H}\psi(h)h^v\big)\right)\\
&=\text{Tr}_{\mathbb F_q(\omega_d)/\mathbb F_q}(e_{\psi})^*={\text{Tr}_{\psi}(e_{\psi})}^*.
\end{align*}
Hence $\mathbb{F}_qG$ is $*$-clean. This completes the proof.
\qed

It is easy to see that Theorem \ref{mainthm} is a generalization of \cite[Theorem 1.1]{HRZ} in the finite field case (the general case will be discussed soon). In \cite[Proof of Theorem 1.1]{HRZ}, there was an inaccuracy about the structure of $Gal(\mathbb F_q(\omega_{n})/\mathbb F_q)$ and we have fixed it in the above proof.

Next, we consider the general case. Firstly, we need to slightly extend the definition of $\sigma_1$. Let $G$ be a finite abelian group with exponent $n$ and $\mathbb{F}$ an arbitrary field of characteristic $p>0$. We define $\sigma_1': \mathbb FG\rightarrow \mathbb FG$ with $\sum_{g\in G}a_gg=\sum_{g\in G}a_gg^v$, where $v^2\equiv 1\pmod{n}$ and $v\not\equiv 1\pmod{n}$. It is easy to see that $Gal(\mathbb{F}(\omega_n)/\mathbb{F})$ is isomorphic to a subgroup of $Gal(\mathbb{F}_p(\omega_n)/\mathbb{F}_p)$; see \cite[Lemma 6.1]{HRZ}. As $Gal(\mathbb{F}_p(\omega_n)/\mathbb{F}_p)$ is generated by the Frobenius map:
$$\Phi:\mathbb{F}_p(\omega_{n})\rightarrow \mathbb{F}_p(\omega_{n})\text{ with }\Phi(x)=x^p,$$
we may assume that $Gal(\mathbb{F}(\omega_n)/\mathbb{F})$ is generated by
$$\Phi:\mathbb{F}(\omega_{n})\rightarrow \mathbb{F}(\omega_{n})\text{ with }\Phi(x)=x^{p^k},$$
for some $k\in\mathbb N$ and let $q=p^k$. Similar to the proof of Theorem \ref{mainthm}, we obtain the following theorem.

\begin{theorem}\label{oldtheorem}
Let $\mathbb F$ be a field of characteristic $p>0$ as described above.
Let $G$ be a finite abelian group with exponent $n$. Let $G=Syl_p(G)\times H$ and $m$ the exponent of $H$, where $Syl_p(G)$ is the Sylow $p$-subgroup of $G$. Taking the involution $\sigma_1'$ as $*$ on $\mathbb{F}G$, then $\mathbb{F}G$ is $*$-clean if and only if there exists $t\in \mathbb N$ such that $q^t\equiv v\pmod{m}.$
\end{theorem}

In the following, we show that there are infinitely many group rings $\mathbb F_q G$ on which the only possible involutions are of the form $\sigma_1$. For simplicity, we consider the case when $G$ is a finite abelian group of odd order. In this case, $\mathbb F_q G$ has at least two different involutions, the identity map and the classical involution. We regard the identity map as an involution of the form $\sigma_1$.

\begin{proposition}\label{all involutions}
Let $\mathbb F_q$ be a finite field of $q$ elements. Let $G$ be a finite abelian group of odd order. Then all involutions of the group ring $\mathbb F_q G$ are of the form $\sigma_1$ if and only if $G$ is a cyclic group $C_p$ of prime order $p$ and one of the following conditions holds:
\begin{enumerate}

\item ord$_p(q)=\frac{p-1}{2}$ and $p\equiv 3\pmod 4$;

\item ord$_p(q)=p-1$.

\end{enumerate}

\end{proposition}

\begin{proof}
Let $p\ge 3$ be a prime. We first prove that all involutions of the group ring $\mathbb F_q C_p$ are of the form $\sigma_1$ if and only if either (1) or (2) holds. Assume that all involutions of the group ring $\mathbb F_q C_p$ are of the form $\sigma_1$. Note that $v^2\equiv 1\pmod{p}$ has only two solutions, $\mathbb F_q C_p$ has only two involutions of the form $\sigma_1$.

Let ord$_p(q)=d$. Then it is easy to see that
$$x^p-1=(x-1)f_1(x)\cdots f_{\frac{p-1}{d}}(x),$$
where $f_i(x)$ $(1\le i\le \frac{p-1}{d})$ is a monic irreducible polynomial of degree $d$ over $\mathbb F_q$. Therefore, we have
$$\mathbb F_q C_p\cong\mathbb F_q[x]/\langle x^p-1\rangle
\cong\mathbb F_{q}\oplus\overbrace{\mathbb F_{q^d}\oplus\cdots\oplus\mathbb F_{q^d}}^{\frac{p-1}{d}}.$$
Note that transpositions between the $\mathbb F_{q^d}$'s are involutions of $\mathbb F_q C_p$. If $\frac{p-1}{d}\ge 3$, then we have at least three transpositions between the $\mathbb F_{q^d}$'s, which means that $\mathbb F_q C_p$ has other forms of involutions and therefore contradicts our assumption. Consequently, we have $\frac{p-1}{d}=1$ or $2$. If $\frac{p-1}{d}=2$, then we must have $d=\frac{p-1}{2}$ is odd, otherwise $\mathbb F_{q^d}\oplus\mathbb F_{q^d}$ has four involutions. Therefore, in this case, we have $p\equiv 3\pmod 4$. If $\frac{p-1}{d}=1$, which means that $\mathbb F_q C_p\cong\mathbb F_{q}\oplus \mathbb F_{q^{p-1}}$. Since $p\ge 3$, it is clear that $\mathbb F_q C_p$ has only two involutions. The converse statement follows from the above discussion.

Now, it suffices to show that if $G$ is not a cyclic group of prime order, then $\mathbb F_q G$ have other forms (except for $\sigma_1$) of involutions. We distinguish several cases.

{\bf Case 1.} $G=C_{p^{m}}$, where $m\ge 2$. Note that $v^2\equiv 1\pmod{p^m}$ has only two solutions, $\mathbb F_q C_{p^m}$ has only two involutions of the form $\sigma_1$. It suffices to prove that $\mathbb F_q G$ has more than two involutions. Let $d_i=\text{ord}_{p^i}(q)$ for $1\le i\le m$.
Similar to the above, it is easy to see that
$$\mathbb F_q C_{p^m}\cong\mathbb F_q[x]/\langle x^{p^m}-1\rangle
\cong\mathbb F_{q}\oplus\overbrace{\mathbb F_{q^{d_1}}\oplus\cdots\oplus\mathbb F_{q^{d_1}}}^{\frac{p-1}{d_1}}\oplus\cdots\oplus
\overbrace{\mathbb F_{q^{d_m}}\oplus\cdots\oplus\mathbb F_{q^{d_m}}}^{\frac{p^m-p^{m-1}}{d_m}}.$$
Also note that, by Lemma \ref{ord_p}, either $\frac{p^m-p^{m-1}}{d_m}=1$ or $p|\frac{p^m-p^{m-1}}{d_m}$. If all involutions of the group ring $\mathbb F_q C_{p^m}$ are of the form $\sigma_1$, then, also by Lemma \ref{ord_p}, we must have $\frac{p-1}{d_1}=\cdots=\frac{p^m-p^{m-1}}{d_m}=1$. In this case, we have
$$\mathbb F_q C_{p^m}
\cong\mathbb F_{q}\oplus\mathbb F_{q^{p-1}}\oplus\cdots\oplus
\mathbb F_{q^{p^m-p^{m-1}}}.$$
As $p-1|\cdots|p^m-p^{m-1}$ are all even and $m\ge 2$, it is clear that $\mathbb F_q C_{p^m}$ has at least four involutions.

{\bf Case 2.} $G=C_{p^{m}}\times C_{p^{n}}$. It suffices to prove that $\mathbb F_q G$ has more than two involutions. Note that
$\mathbb F_q G=[\mathbb F_q C_{p^{m}}]C_{p^{n}}$. If $m\ge 2$, then $\mathbb F_q C_{p^{m}}\cong \mathbb K_1\oplus\cdots\oplus \mathbb K_t$, where $t\ge 3$ and each $\mathbb K_i$ is a finite extension of $\mathbb F_q$. Therefore
$$\mathbb F_q G\cong[\mathbb K_1\oplus\cdots\oplus \mathbb K_t]C_{p^{n}}
\cong\mathbb K_1C_{p^{n}}\oplus\cdots\oplus \mathbb K_tC_{p^{n}}.$$
Consequently, $\mathbb F_q G$ has at least 8 involutions. It follows that $m=n=1$. Based on our discussion above, it suffices to consider that $\mathbb F_q C_p\cong\mathbb F_{q}\oplus \mathbb F_{q^{p-1}}$ and
$$\mathbb F_q G\cong[\mathbb F_{q}\oplus \mathbb F_{q^{p-1}}]C_p
\cong\mathbb F_{q}C_p\oplus \mathbb F_{q^{p-1}}C_p.$$
As $p-1$ is even, it is clear that $\mathbb F_{q^{p-1}}C_p$ (hence $\mathbb F_q G$) has at least four involutions. The above discussion can be modified to show that if $G$ is a finite abelian $p$-group and $G\neq C_p$, then $\mathbb F_q G$ has at least three involutions.

{\bf Case 3.} $G=H_{p_1}\times H_{p_2}$, where $H_{p_1}$ (resp. $H_{p_1}$) is a finite abelian $p_1$-group (resp. $p_2$-group) with $p_1\neq p_2$. Then $\mathbb F_q G=[\mathbb F_q H_{p_1}]H_{p_2}$. Note that, in this case, $v^2\equiv 1\pmod{\exp(G)}$ has only four solutions, $\mathbb F_q G$ has only four involutions of the form $\sigma_1$. It suffices to prove that $\mathbb F_q G$ has more than four involutions. Based on our discussion above, we have either
$\mathbb F_q H_{p_1}\cong\mathbb F_{q}\oplus \mathbb F_{q^{p_1-1}}$ or
$\mathbb F_q H_{p_1}=\mathbb K_1\oplus\cdots\oplus \mathbb K_t$, where $t\ge 3$ and each $\mathbb K_i$ is a finite extension of $\mathbb F_q$. In either case, it is easy to see that $\mathbb F_q G=[\mathbb F_q H_{p_1}]H_{p_2}$ has more than four involutions, as desired.

Since any finite abelian group $G$ can be decomposed as the direct product of $p$-groups. The desired result can be easily verified by modifying the above approach.
\end{proof}

We also provide an elementary method to study the $*$-cleanness of $\mathbb F_2G$ as well as some examples.

We first show that the idempotents in $\mathbb F_2G$ can be easily described. Actually, let $a=g_1+\cdots+g_k\in\mathbb F_2G$, with $g_i\in G$. Then, $a$ is an idempotent if
$$(g_1+\cdots+g_k)^2=g_1^2+\cdots+g_k^2=g_1+\cdots+g_k.$$
This means that $\{g_1^2,\ldots,g_k^2\}=\{g_1,\ldots,g_k\}$. Moreover, if $g\in \{g_1,\ldots,g_k\}$, then so is $g^{2^i}$ for any $i\in\mathbb N$. Consequently, we can write $a$ in the following form
\begin{equation}\label{idemform}
a=(g_{i_1}+g_{i_1}^2+\cdots+g_{i_1}^{2^{l_1}})+\cdots
+(g_{i_s}+g_{i_s}^2+\cdots+g_{i_s}^{2^{l_s}}),
\end{equation}
where the $g_{i_j}$'s are distinct and the $l_j$'s are the smallest integers such that $g_{i_j}^{2^{l_j+1}}=g_{i_j}$.

We consider the classical involution $*$ on $\mathbb F_2G$ and Theorem \ref{mainthm} is reduced to the following one.

\begin{theorem}\label{simpleF2}
Let $G$ be a finite abelian group with exponent $n$ and $(2,n)=1$. Then $\mathbb{F}_2G$ is $*$-clean if and only if there exists $t\in\mathbb N$ such that $2^t\equiv -1\pmod{n}$.
\end{theorem}

\begin{proof}
Firstly, we assume that $\mathbb F_2G$ is $*$-clean. Therefore, every idempotent is a projection. In particular, let $g\in G$ with ord$(g)=n$, we consider the following idempotent:
$$a=g+g^2+\cdots+g^{2^k},$$
where $k$ is the smallest positive integer such that $g^{2^{k+1}}=g$. Since $a^*=g^{-1}+g^{-2}+\cdots+g^{-2^k}$ and $a^*=a$, we have $g^{-1}\in \{g,g^2,\ldots,g^{2^k}\}$. Therefore, there exists a positive integer $t$ such that $2^t\equiv -1\pmod{n}$.

Conversely, we assume that there exists $t\in\mathbb N$ such that $2^t\equiv -1\pmod{n}$. By Lemmas \ref{lemma3.1} and \ref{cleanFG}, it suffices to prove that any idempotent $a$ is a projection. Based on the above discussion, any idempotent can be written in the form (\ref{idemform}). Therefore, we have
$$a^*=(g_1^{-1}+g_1^{-2}+\cdots+g_1^{-2^{l_1}})+\cdots
+(g_s^{-1}+g_s^{-2}+\cdots+g_s^{-2^{l_s}}).$$
As there exists $t\in\mathbb N$ such that $2^t\equiv -1\pmod{n}$, it is easy to see that $a^*=a$. This completes the proof.
\end{proof}

Note that, the above approach is also valid if we consider involutions of the form $\sigma_1$ on $\mathbb F_2G$ instead. Next, we provide some examples.

\begin{example}
Consider $\mathbb F_2(C_3\times C_9)$. As $2^3\equiv -1\pmod{9}$,
by Theorem \ref{simpleF2}, $\mathbb F_2(C_3\times C_9)$ is $*$-clean.
\end{example}

\begin{example}
Consider $\mathbb F_2(C_3\times C_{15})$. It is easy to verify that there exists no $t\in\mathbb N$ such that $2^t\equiv -1\pmod{15}$. Therefore, by Theorem \ref{simpleF2}, $\mathbb F_2(C_3\times C_{15})$ is not $*$-clean.
\end{example}

\medskip

Now we turn to Theorem \ref{mainthm2}. Since the proof of Theorem \ref{mainthm2} is very similar to the above proof of Theorem \ref{mainthm}, we omit it here. We provide two examples when the involution $*$ is taken as: $\mathbb F_{q^2}G\rightarrow \mathbb F_{q^2}G$ with $\sum_{g\in G}a_gg\mapsto\sum_{g\in G}a_g^qg^{-1}$.

\begin{example}
Consider $\mathbb F_{2^2}(C_9\times C_{9})$. As $2^{2\cdot2}\equiv -2\pmod{9}$, by Theorem \ref{mainthm2}, $\mathbb F_{2^2}(C_9\times C_{9})$ is $*$-clean.
\end{example}

\begin{example}
Consider $\mathbb F_{2^2}(C_5\times C_{25})$. It is easy to verify that there exists no $t\in\mathbb N$ such that $2^{2t}\equiv -2\pmod{25}$. Therefore, by Theorem \ref{mainthm2}, $\mathbb F_{2^2}(C_5\times C_{25})$ is not $*$-clean.
\end{example}

\section{The connection between the $*$-cleanness of $\mathbb F_qG$ and coding theory}

In this section, based on the construction at the beginning of Section 3 and the proof of Theorem \ref{mainthm}, we provide a brief discussion on the connection between the $*$-cleanness (with the classical involution $*$) of $\mathbb F_qG$ and abelian group codes in $\mathbb F_qG$. Here, $G$ is a finite abelian group with exponent $n$ and $(q,n)=1$. By Lemma \ref{PIC}, we define an equivalence relation on $\hat{G}$: $\psi_1\sim\psi_2$ if and only if they belong to the same $q$-cyclotomic class in $\hat{G}$. Let $\overline{\psi}$ be the equivalent class of $\psi$ and $\overline{\hat{G}}$ a set of representatives of $\hat{G}/\sim$.

\begin{lemma}\label{CharCbot}
Let $\mathcal C_{\psi}$ be the abelian code generated by $\text{Tr}_{\psi}(e_\psi)$ in $\mathbb F_qG$. Then we have
$$\mathcal C_{\psi}^{\bot}=\bigoplus_{\varphi\in\overline{\hat{G}}
\setminus\{\overline{\psi^{-1}}\}}\mathcal C_{\varphi},$$
where $\mathcal C_{\varphi}$ is the abelian code generated by $\text{Tr}_{\varphi}(e_\varphi)$ in $\mathbb F_qG$.
\end{lemma}

\begin{proof}
Let $\alpha=\sum_{g\in G}a_g g$ and $\beta=\sum_{g\in G}b_g g$ in $\mathbb F_qG$. Then $\langle\alpha,\beta\rangle=0$ is equivalent to $\langle\alpha\beta^*,1\rangle=0$, where $\beta^*=\sum_{g\in G}b_g g^{-1}$. Note that, $\text{Tr}_{\varphi}(e_\varphi)^*=\text{Tr}_{\varphi^{-1}}(e_{\varphi^{-1}})$. Therefore, it is easy to see that
$\langle\text{Tr}_{\varphi}(e_\varphi),\text{Tr}_{\psi}(e_\psi)\rangle=0$ for all $\varphi\in\overline{\hat{G}}\setminus\{\overline{\psi^{-1}}\}$. It follows that
$\bigoplus_{\varphi\in\overline{\hat{G}}
\setminus\{\overline{\psi^{-1}}\}}\mathcal C_{\varphi}\subset
\mathcal C_{\psi}^{\bot}$. If $\bigoplus_{\varphi\in\overline{\hat{G}}
\setminus\{\overline{\psi^{-1}}\}}\mathcal C_{\varphi}\subsetneq
\mathcal C_{\psi}^{\bot}$, then we have $\mathcal C_{\psi}^{\bot}=\bigoplus_{\varphi\in\overline{\hat{G}}}\mathcal C_{\varphi}=\mathbb F_qG$ which leads to $\mathcal C_{\psi}=0$, a contradiction.
\end{proof}

From the proof of Theorem \ref{mainthm}, we have the following results.

\begin{proposition}\label{LCDortho}
Let $\psi\in\hat{G}$ and ord$(\psi)=d$. Let $\mathcal C_{\psi}$ be the abelian code generated by $\text{Tr}_{\psi}(e_\psi)$ in $\mathbb F_qG$. Then we have:
\begin{enumerate}

\item $\mathcal C_{\psi}$ is an LCD abelian code if and only if $\psi\sim\psi^{-1}$, or equivalently, there exists $t\in \mathbb N$ such that $q^t\equiv -1\pmod{d}$;
\item $\mathcal C_{\psi}$ is a self-orthogonal code if and only if $\psi\nsim\psi^{-1}$, or equivalently, there exists no $t\in\mathbb N$ such that $q^t\equiv -1\pmod d$.
\end{enumerate}
\end{proposition}

Combining Theorem \ref{mainthm} and Proposition \ref{LCDortho}, we have the following theorem.

\begin{theorem}\label{starcodes}
The group ring $\mathbb F_qG$ is $*$-clean if and only if one of the following holds:
\begin{enumerate}
\item all abelian group codes are LCD abelian codes in $\mathbb F_qG$;
\item all abelian group codes are not self-orthogonal abelian codes in $\mathbb F_qG$;
\item if there exists an LCD abelian codes $\mathcal C_{\psi}$ with ord$(\psi)=n$ in $\mathbb F_qG$;
\item if there exists $\psi\in\hat{G}$ with ord$(\psi)=n$ such that $\mathcal C_{\psi}$ is not a self-orthogonal abelian code in $\mathbb F_qG$.
\end{enumerate}
\end{theorem}

Finally, we provide some examples.

\begin{example}
Consider $\mathbb F_2(C_3\times C_3)$ and let $G=C_3\times C_3=\langle x\rangle\times\langle y\rangle$. Let $\psi\in\hat{G}$ which is defined by $\psi(x)=1$ and $\psi(y)=\omega$, where $\omega$ is a third primitive root of unity over $\mathbb F_2$.
Let $\mathcal C_{\psi}$ be the abelian code generated by $\text{Tr}_{\psi}(e_\psi)$. On the one hand, as $2\equiv -1\pmod 3$, it follows immediately from Proposition \ref{LCDortho}.(1) that $\mathcal C_{\psi}$ is an LCD abelian code in $\mathbb F_2G$. On the other hand, a direct computation shows that $$\text{Tr}_{\psi}(e_\psi)=(1+x+x^2)(y+y^2)$$ and
$$\mathcal C_{\psi}=\{0,\ \text{Tr}_{\psi}(e_\psi),\ \text{Tr}_{\psi}(e_\psi)y,\ \text{Tr}_{\psi}(e_\psi)y^2\}.$$
It is easy to see that for $i\neq j$, we have
\begin{align*}
&\langle\text{Tr}_{\psi}(e_\psi)y^i,\text{Tr}_{\psi}(e_\psi)y^j\rangle\\
&=\langle(1+x+x^2)(y^{1+i}+y^{2+i}),(1+x+x^2)(y^{1+j}+y^{2+j})\rangle=1.
\end{align*}
Consequently, we have
$\mathcal C_{\psi}\cap\mathcal C_{\psi}^{\bot}=\{0\}$ and $\mathcal C_{\psi}$ is an LCD abelian code in $\mathbb F_2G$.
\end{example}

\begin{example}
Consider $\mathbb F_2(C_7\times C_7)$ and let $G=C_7\times C_7=\langle x\rangle\times\langle y\rangle$. Let $\psi\in\hat{G}$ which is defined by $\psi(x)=1$ and $\psi(y)=\omega$, where $\omega$ is a 7-th primitive root of unity over $\mathbb F_2$ and its minimal polynomial over $\mathbb F_2$ is $x^3+x+1$. Let $\mathcal C_{\psi}$ be the abelian code generated by $\text{Tr}_{\psi}(e_\psi)$. On the one hand, as there exists no $t\in\mathbb N$ such that $2^t\equiv -1\pmod 7$, it follows immediately from Proposition \ref{LCDortho}.(2) that $\mathcal C_{\psi}$ is a self-orthogonal abelian code in $\mathbb F_2G$. On the other hand, a direct computation shows that
$$\text{Tr}_{\psi}(e_\psi)=(\sum_{i=0}^6x^i)(1+y^3+y^5+y^6),$$
and
$$\mathcal C_{\psi}=\{0,\ \text{Tr}_{\psi}(e_\psi),\ \text{Tr}_{\psi}(e_\psi)y,
\ldots,\text{Tr}_{\psi}(e_\psi)y^6\}.$$
It is easy to see that for $i\neq j$, we have
$$|\{\overline{i},\ \overline{3+i},\ \overline{5+i},\ \overline{6+i}\}
\cap\{\overline{j},\ \overline{3+j},\ \overline{5+j},\ \overline{6+j}\}|=2,$$
where $\overline{a}\equiv a\pmod 7$ for any $a\in\mathbb N$. Therefore, we have
\begin{align*}
&\langle\text{Tr}_{\psi}(e_\psi)y^i,\text{Tr}_{\psi}(e_\psi)y^j\rangle\\
&=\langle(\sum_{i=0}^6x^i)(y^{i}+y^{3+i}+y^{5+i}+y^{6+i}),
(\sum_{i=0}^6x^i)(y^{j}+y^{3+j}+y^{5+j}+y^{6+j})\rangle=0.
\end{align*}
Consequently, we have $\mathcal C_{\psi}\subset\mathcal C_{\psi}^{\bot}$ and therefore $\mathcal C_{\psi}$ is a self-orthogonal abelian code in $\mathbb F_2G$.
\end{example}

\subsection*{Acknowledgments}

We sincerely thank the referees for their very helpful comments. In particular, Proposition \ref{all involutions} was inspired by an insightful comment from a referee. D.C. Han was supported by the National Science Foundation of China Grant No.11601448 and the Fundamental Research Funds for the Central Universities Grant No.2682020ZT101. H.B. Zhang was supported by the National Science Foundation of China Grant No.11901563.

\end{document}